\documentclass[journal]{IEEEtran}
\usepackage{amssymb,amsthm,bm}
\usepackage{mathtools}
\usepackage{graphicx}
\usepackage{booktabs,multirow,array,cite}
\usepackage{algorithm,algpseudocode}
\usepackage{siunitx,url,color,enumitem}
\usepackage{tikz,physics}
\usepackage{microtype}   
\usetikzlibrary{positioning,arrows.meta,calc,fit}
\tikzset{>=Latex}
\allowdisplaybreaks

\newtheorem{assumption}{Assumption}
\newtheorem{proposition}{Proposition}

\DeclareMathOperator{\Range}{Range}
\DeclareMathOperator{\diag}{diag}

\makeatletter
\@ifundefined{norm}{%
	\DeclarePairedDelimiter{\norm}{\lVert}{\rVert}%
}{%
	\let\norm\mynorm
}
\makeatother

\begin{document}
	
 \title{Co-State Based Data Fusion and Risk Aware Filtering for Spacecraft Navigation and Hazard Prediction}
 \author{Surya Ratna Prakash D \thanks{U R Rao Satellite Centre, Indian Space Research Organisation, Bengaluru 560017, India and Computational and Data Sciences,  Indian Institute of Science, Bengaluru 560012, India. Email: dsrp@ursc.gov.in} \and  Soumyendu Raha
\thanks{Computational and Data Sciences, Indian Institute of Science, Bengaluru 560012. Email: raha@iisc.ac.in}}

\maketitle

\begin{abstract}
	This paper develops a co-state based fusion framework for
	spacecraft navigation, consistency monitoring, and hazard forecasting. A
	differential--algebraic co-state is introduced as an instantaneous Lagrange
	multiplier that enforces measurement--dynamics compatibility at the differential
	level and provides a physically interpretable signal of geometric inconsistency.
	On a longer time scale, co-state and innovation trajectories are used to learn a
	continuous-time Markov generator governing transitions between coarse behavioral
	regimes, enabling intrinsic probabilistic risk forecasting through mode
	probabilities and mean first-passage time (MFPT).
	
	The resulting architecture unifies geometric projection, stochastic inference,
	and probabilistic risk assessment in a single online pipeline without requiring
	predefined fault models, labelled failure data, or heuristic thresholds. The
	framework is demonstrated on real lunar powered-descent telemetry, 
	where it detects structural internal model inconsistency significantly earlier than physical
	divergence or statistical inconsistency in an Extended Kalman Filter (EKF).
	The results show that geometric inconsistency, stochastic drift, and probabilistic
	risk rise coherently prior to failure, yielding interpretable and operationally
	meaningful early-warning capability for autonomous landing systems.
\end{abstract}

\begin{IEEEkeywords}
	Spacecraft navigation, co-state dynamics, stochastic generator matrix, nonlinear filtering, hazard prediction, risk-aware estimation, differential--algebraic systems.
\end{IEEEkeywords}

\section{Introduction}

Autonomous spacecraft must operate safely under increasingly complex dynamics,
uncertain environments, and limited onboard sensing and computation. This is particularly critical than in powered descent and landing, where small modeling
errors, sensor inconsistencies, or guidance mismatches can lead to catastrophic
outcomes within short time horizons. Ensuring not only accurate state estimation,
but also early detection of internal inconsistency and impending failure, is
therefore a central challenge for autonomous navigation and control.

Classical navigation architectures rely primarily on stochastic state estimation
methods such as the Kalman filter and its nonlinear variants. These filters are
designed to optimally fuse noisy measurements and dynamics models under
probabilistic assumptions, and they remain the backbone of modern aerospace
navigation systems. However, they enforce consistency between measurements and
models only in a statistical sense, through innovation correction and covariance
propagation. When the assumed model becomes internally inconsistent with the
physical system or sensing geometry, the filter may remain numerically stable
while absorbing the mismatch into its uncertainty representation, thereby
masking early indicators of failure.

At the same time, fault detection, anomaly monitoring, and risk-aware autonomy
have become increasingly important in aerospace systems. These approaches
typically rely on residual-based detectors, mode-switching logic, or learned
classifiers, and often require explicit fault models, heuristic thresholds, or
labeled data. While effective in specific contexts, such methods do not directly
expose the geometric structure of measurement--dynamics compatibility, nor do
they naturally provide intrinsic long-horizon risk measures tied to the
underlying system geometry.

This paper introduces a co-state based fusion framework that
explicitly enforces measurement--dynamics consistency at the differential level
and simultaneously captures slow-time hazard evolution through data-driven
stochastic modeling. The central idea is to introduce an algebraic co-state
variable that acts as an instantaneous Lagrange multiplier, projecting the state
evolution onto the measurement-consistent tangent manifold. The magnitude of this
co-state directly quantifies instantaneous geometric incompatibility between the
assumed dynamics and observed measurements, providing a physically interpretable
consistency signal that is absent from classical filters.

On a longer time scale, the co-state and innovation signals are used to learn a
continuous-time Markov generator governing transitions between coarse behavioral
regimes. This generator induces a probabilistic semigroup over modes, enabling
computation of intrinsic risk measures such as hazard probability and mean
first-passage time (MFPT) without requiring predefined fault models, failure
labels, or heuristic thresholds.

The resulting framework unifies geometric projection, stochastic inference, and
probabilistic hazard forecasting into a single operational pipeline that is
suitable for online, asynchronous, and embedded implementation. It does not
replace classical state estimation, but rather augments it with explicit geometric
consistency monitoring and predictive risk awareness.

The main contributions of this work are:
\begin{itemize}
	\item A deterministic co-state differential--algebraic formulation that enforces
	instantaneous measurement--dynamics compatibility and yields a geometric
	inconsistency signal with clear physical interpretation.
	\item A data-driven stochastic drift--diffusion and generator learning framework
	that captures slow-time regime evolution and enables intrinsic probabilistic risk
	forecasting via mode probabilities and mean first-passage time.
	\item A coupled deterministic--stochastic fusion pipeline that produces early
	risk-aware warning signals without labeled data or predefined fault hypotheses.
	\item A demonstration on a lunar lander's powered-descent telemetry showing that
	internal model inconsistency is detected significantly earlier than physical
	divergence or Extended Kalman Filter (EKF) based statistical alarms.
\end{itemize}

The remainder of the paper is organized as follows. Section~\ref{sec:related}
reviews related work. Sections~\ref{sec:sdamodel} and~\ref{sec:approx-meas}
introduce the deterministic co-state formulation and measurement model.
Sections~\ref{sec:Drift--Diffusion} and~\ref{sec:AdaptiveWeighting} develop the
stochastic drift--diffusion model and generator learning. Section~\ref{sec:Algorithm}
presents the algorithmic implementation. Section~\ref{sec:results-lander}
demonstrates the method on lunar descent telemetry and analyzes the results. Finally, Section~\ref{sec:DiscussionFuturework} concludes with discussion and
future directions.

\section{Related Work}
\label{sec:related}

The problem of safe and autonomous spacecraft descent lies at the intersection of stochastic state estimation, fault detection, and probabilistic risk-aware decision-making. The proposed framework draws on and extends ideas from each of these areas.

\subsection{Model-Based Estimation and Filtering}

Spacecraft navigation has traditionally relied on Bayesian filtering methods such as the Kalman filter, EKF, Unscented Kalman Filter, and particle filters \cite{kalman1960,kalmanbucy1961,doucet2001}. These methods enforce measurement consistency statistically through innovation correction and covariance propagation.

While highly effective for state estimation, these filters treat measurement--dynamics consistency implicitly and in expectation. Model mismatch, partial observability, or ill-conditioned sensing geometry are absorbed into the innovation and covariance, often obscuring early indicators of internal inconsistency or impending failure.

\subsection{Consistency Monitoring and Fault Detection}

Fault detection in aerospace systems has traditionally relied on residual monitoring, hypothesis testing, parity relations, and change detection techniques. More recently, data-driven anomaly detection and machine learning approaches have been explored.

These methods often require labeled failure data, offline training, or heuristic thresholds and may lack physical interpretability. In contrast, the co-state formulation introduced here provides a physically interpretable consistency signal that directly measures geometric incompatibility between dynamics and measurements at the differential level.

\subsection{Risk-Aware Navigation and Autonomous Landing}

Recent work emphasizes robustness and risk-aware decision-making for autonomous landing
\cite{GaoZhou2018_TerrainHazard,Calkins2022_RobustTrajectoryPDL,JungBang2021_FaultTolerantMPC,Fourlas2021_FaultDiagnosisFTC,WangPengJiang2020_RT_FaultDetection}.
These approaches typically rely on residual-based detectors, interacting multiple-model filters, adaptive covariance tuning, or learned transition models.

However, consistency between measurements and dynamics is generally handled implicitly through statistical correction, mode switching, or threshold logic rather than being enforced explicitly as a geometric constraint on state evolution.

\subsection{Stochastic Hybrid Systems and Mode-Based Risk}

Markov jump systems, hidden Markov models, and stochastic hybrid systems enable probabilistic reasoning over discrete modes and hazard states. These models typically assume predefined modes or fault hypotheses and require labeled data.

In contrast, the present framework learns a continuous-time generator directly from observed geometric and innovation features without requiring predefined fault models or failure labels.

\subsection{Positioning of the Present Work}

The proposed framework differs from existing approaches in three key respects:
\begin{itemize}
	\item It introduces an explicit algebraic co-state enforcing measurement--dynamics compatibility at the differential level.
	\item It couples this geometric signal to a learned continuous-time generator enabling intrinsic risk forecasting.
	\item It integrates geometric, stochastic, and probabilistic components into a unified online architecture.
\end{itemize}

Accordingly, the present work should be viewed not as a replacement for classical filtering or hazard detection methods, but as a complementary layer that exposes latent internal inconsistency and converts it into interpretable and predictive risk information.

\section{Stochastic Differential--Algebraic Equation Model}
\label{sec:sdamodel}
Let the spacecraft state $x_t \in \mathbb{R}^n$ evolve according to
\begin{equation}
	\dot{x}_t = f(x_t,t), \label{eq:dyn}
\end{equation}
which, once observation noise is introduced, is interpreted in the It\^{o} differential sense. The $m$-dimensional measurement model is
\begin{equation}
	y_t^{\mathrm{model}} := h(x_t,t) + w_t, \quad w_t \sim \mathcal{N}(0,R_t). \label{eq:meas}
\end{equation}
Differentiating \eqref{eq:meas} yields
\begin{equation}
	dy_t^{\mathrm{model}} := H_t\,dx_t + \bigg(\frac{\partial h}{\partial t}\bigg)_t dt + dw_t,
	\quad H_t := \bigg(\frac{\partial h}{\partial x}\bigg)_t. \label{eq:dmeas}
\end{equation}
where $f:{\mathbb R}^n \times {\mathbb R}_{+}\to {\mathbb R}^n$ describes the system dynamics, $h:{\mathbb R}^n \times {\mathbb R}_{+} \to {\mathbb R}^m$ is the measurement function, and $y^{\mathrm{observed}}_t, y^{\mathrm{model}}_t$ for $t \ge 0$ are the actual observation and the model observation processes respectively. 

Measurement--dynamics consistency requires that discrepancies between the observed and model-implied measurement increments be explainable by admissible state variations, i.e.,
\begin{equation}
	H_t \big(dx_t - f_t dt\big) = dy_t^{\mathrm{observed}} - dy_t^{\mathrm{model}}. \label{eq:cons}
\end{equation}
Equivalently, projecting onto the range of $H_t^\top$ (the measurement tangent space),
\begin{equation}
	dx_t - f_t dt - H_t^{\mathsf T}(H_t H_t^{\mathsf T})^{-1}
	\big(dy_t^{\mathrm{observed}} - dy_t^{\mathrm{model}}\big) = 0. \label{eq:measmodel}
\end{equation}

Introduce an algebraic co-state $\lambda_t \in \mathbb{R}^m$ acting as an instantaneous Lagrange multiplier, and augment the dynamics as
\begin{subequations}\label{eq:sdae}
	\begin{align}
		dx_t &= f_t dt + H_t^{\mathsf T} \lambda_t dt, \label{eq:aug}\\
		dy_t^{\mathrm{observed}} &= \Big(H_t f_t + \big(\tfrac{\partial h}{\partial t}\big)_t\Big)dt
		+ (H_t H_t^{\mathsf T})\lambda_t dt + dw_t. \label{eq:yproj}
	\end{align}
\end{subequations}
Under regularization $\varepsilon > 0$, 
this system can be interpreted as an index-1 stochastic differential--algebraic equation.

Thus the algebraic co-state satisfies, in the sense of differential--algebraic consistency,
\begin{equation}
	(H_t H_t^{\mathsf T})\lambda_t dt
	= dy_t^{\mathrm{observed}} - \eta_t dt - dw_t, 
	\quad \eta_t :=\big(\frac{\partial h}{\partial t}\big)_t + H_t f_t. \label{eq:lamraw}
\end{equation}

When $H_tH_t^{\mathsf T}$ is ill-conditioned, we regularize and include whitening:
\begin{equation}
	\lambda_t dt = (H_tH_t^{\mathsf T} + \varepsilon I_m)^{-1}\Sigma_t^{-1}dy_t^{\mathrm{observed}} - \eta_t dt - \sigma_t dW_t, \label{eq:lamreg}
\end{equation}
with $\varepsilon>0$ spectral regularizer, $\Sigma_t$ the innovation weighting (whitening) and $\sigma_t$ the observation noise shaping matrix (see Sec.~\ref{sec:wtobs}). The ratio is interpreted in the incremental sense over finite sampling intervals $\Delta t$, rather than as a pointwise derivative, consistent with the discrete-time realization of the continuous-time model.

\section{Weighted Observations}
\label{sec:wtobs}

To improve numerical conditioning and enable interpretable innovation statistics, we introduce a weighted $m$-dimensional observation model by shaping the Wiener process in \eqref{eq:dmeas} with a full-rank matrix $\sigma_t = \sigma(x_t,t) \in \mathbb{R}^{m \times m}$, which may depend on the state $x_t$ and time $t$. Importantly, this weighting does not alter the physical observation process; rather, it defines a surrogate weighted observation model used solely for numerical conditioning and statistical normalization.

Under this weighted representation, the observation increment is written as
\begin{equation}
	dy_t^{\mathrm{weighted}} := \sigma_t \, dW_t + \eta_t \, dt,
	\label{eq:wtobs}
\end{equation}
where $\eta_t := \big(\frac{\partial h}{\partial t}\big)_t  + H_t f_t$ is the model-predicted measurement increment.

All quantities are interpreted causally with respect to the observed measurement history. Under this interpretation, the algebraic co-state correction assumes the regularized and whitened form used in \eqref{eq:lamreg}, with $\Sigma_t := \sigma_t \sigma_t^{\mathsf T}$ serving as an innovation weighting (whitening) matrix.

This weighting plays a critical operational role in Algorithm~\ref{alg:costate} by:
(i) stabilizing the co-state computation under ill-conditioned sensing geometry,
(ii) producing approximately whitened innovation statistics for regime clustering, and
(iii) enabling consistent downstream generator learning and risk assessment.

\section{Approximate Design of the Measurement Model}
\label{sec:approx-meas}

Accurate measurement modeling is central to deterministic--stochastic
co-state fusion, since the co-state multipliers are defined by the
instantaneous geometric discrepancy between observed measurement
increments and model-implied increments. This section develops a
framework for designing tractable measurement models under:
(i) nonlinear geometric sensing,
(ii) partial and rank-deficient observability,
(iii) nonuniform sampling, and
(iv) stochastic small-noise approximations.
We then derive a consistent approximate model specialized to the lunar
lander sensing suite used in our experiments.

\subsection{General Measurement Map and Differential Approximation}

Let the true measurement process be generated by an unknown nonlinear map
\[
y_t^{\mathrm{observed}} = h_{\rm true}(x_t,t) + \nu_t,
\]
where $\nu_t$ denotes sensor noise, latency effects, and unmodeled bias.
In practice, one replaces $h_{\rm true}$ with a tractable approximation
$h$, yielding the computable model \eqref{eq:meas}
with an $R_t$ that approximates the local covariance of the sensor noise
process.

The co-state formulation operates on the \emph{differential} of
\eqref{eq:meas}, which under a first-order expansion yields \eqref{eq:dmeas}.

This differential approximation assumes:
(i) local smoothness of $h$,
(ii) small noise $dw_t$ with covariance $R_t\,dt$, and
(iii) bounded curvature so that the second-order remainder satisfies
$\|r_t\| = \mathcal O(\|dx_t\|^2)$.

\subsection{Approximate Innovation Modeling}
The drift of the model-predicted measurement
increment is $\eta_t\,dt$, while the observed
increment is $dy_t^{\mathrm{observed}} := y_{t+dt}-y_t$. The innovation is

\begin{equation}
	\tilde{v}_t
	= dy_t^{\mathrm{observed}} - dy_t^{\mathrm{model}}
	= H_t (dx_t - f_t dt) + dw_t - \underbrace{r_t}_{\text{Taylor remainder}},
\end{equation}

where the second-order remainder
\(
r_t = \tfrac12 d x_t^\top
\frac{\partial^2 h}{\partial x^2}
(x_t,t)
d x_t
\)
is non-negligible for large increments, rapid motion, or
strongly nonlinear sensors.

To preserve computational tractability and enable instantaneous
projection, we adopt the small-noise operating regime in which
$r_t \approx 0$ and treat the innovation incrementally as
\begin{equation}
	\tilde{v}_t \approx dw_t
	\quad\Longrightarrow\quad
	\tilde{v}_t \sim \mathcal N(0, R_t\,dt).
\end{equation}
This approximation defines the regime in which Algorithm~\ref{alg:costate} operates,
leads to the geometric consistency condition, and is supported empirically in Section~\ref{sec:results-lander}. 

\begin{equation}
	dx_t - f_t dt \in \Range(H_t^{\mathsf T}).
\end{equation}

\subsection{Rank-Robust Approximate Jacobian}

In practical settings, $H_t$ may be ill-conditioned or rank-deficient
(e.g., near nadir pointing or loss of geometric sensitivity).
We employ the rank-robust approximation
\begin{equation}
	(H_t H_t^\top)^{-1}
	\approx (H_t H_t^\top + \varepsilon I_m)^{-1},
\end{equation}
with $\varepsilon>0$ chosen by adaptive spectral regularization.
This construction approximates the Moore--Penrose projector onto
$\Range(H_t^{\mathsf T})$ while ensuring numerical stability, bounded
co-state updates, and well-defined projection under partial observability.

\subsection{Approximate Noise Shaping and Whitening}

We introduce a surrogate innovation weighting
$\Sigma_t = \diag(\sigma_{t,1}^2,\ldots,\sigma_{t,m}^2)$ that approximates
the instantaneous covariance of $dy_t^{\mathrm{observed}}$.
A robust choice is the rolling RMS of the innovation,
\begin{equation}
	\sigma_{t,i}^2 \approx
	\max\!\Big(\sigma_{\min}^2,\,
	\mathrm{RMS}\{\tilde v_{s,i}:\, s\in[t-\Delta,t]\}^2\Big).
\end{equation}
The resulting whitened innovation,
\begin{equation}
	z_t :=
	\sqrt{
		(dy_t^{\mathrm{observed}} - \eta_t dt)^\top
		\Sigma_t^{-1}
		(dy_t^{\mathrm{observed}} - \eta_t dt)
	},
\end{equation}
is approximately $\chi_m$-distributed under the small-noise assumption
and is used for regime clustering and risk assessment in Algorithm~\ref{alg:costate}.

\subsection{Specialization to the Lunar Lander Measurement Suite}

For the lunar lander dataset, we adopt the approximate measurement model
\begin{equation}
	h(x)
	=
	\begin{bmatrix}
		z \\ r(x) \\ v_z
	\end{bmatrix},
	\qquad
	r(x) := \|p\|_2 = \sqrt{x^2+y^2+z^2},
	\label{eq:lander-meas}
\end{equation}
with Jacobian
\begin{equation}
	H(x)
	=
	\begin{bmatrix}
		0 & 0 & 1 & 0 & 0 & 0\\
		x/r & y/r & z/r & 0 & 0 & 0\\
		0 & 0 & 0 & 0 & 0 & 1
	\end{bmatrix},
	\label{eq:lander-H}
\end{equation}
where $r$ is replaced by $\max(r,\epsilon_r)$ for numerical stability.

\paragraph*{Predicted measurement derivative.}
With kinematic prior $\dot p=v,\ \dot v=0$,
\begin{equation}
	\eta_t
	=
	\begin{bmatrix}
		v_z \\
		\dfrac{x v_x + y v_y + z v_z}{r} \\
		0
	\end{bmatrix}.
\end{equation}

\paragraph*{Approximate co-state law.}
The projected co-state law specialized to this measurement set is
\begin{equation}
	\lambda_t ~dt
	\approx
	(H_t H_t^\top + \varepsilon I)^{-1}
	\Sigma_t^{-1}
	{dy_t^{\mathrm{observed}} - \eta_t dt},
\end{equation}
which is well defined except at $p=0$, where $\varepsilon$ regularizes the
inverse.


\section{Drift, Diffusion and Generator Learning}
\label{sec:Drift--Diffusion}
Let samples $\{(x_i,\lambda_i,t_i)\}$ be clustered into modes
$\{C_k\}_{k=1}^K$, defining a coarse-grained representation of the
underlying continuous dynamics. For transitions $C_k\to C_l$ with index
set $T_{kl}$, define the empirical drift intensity
\begin{align}
	a_{kl}^{(q)} &=
	\frac{1}{|T_{kl}|}
	\sum_{i\in T_{kl}}
	\frac{\big|x_{i+1}^{(q)}-x_i^{(q)}\big|}{t_{i+1}-t_i}, \\
	\bar{a}_k^{(q)} &=
	\sum_{l}
	\frac{a_{kl}^{(q)}}{|l_{kl}^{(q)}|},
\end{align}
where absolute increments are used to estimate transition intensities between modes, independent of signed direction within a mode.

Aggregate diffusion statistics $\bar{\Sigma}_k$ and assemble the
continuous-time generator $L$ as
\begin{align}
	L_{kl}=\bar{a}_{kl}\ (k\neq l), \quad
	L_{kk}=-\sum_{j\neq k}L_{kj}, \quad
	\mathbf{1}^{\mathsf T}L=0.
	\label{eq:gen}
\end{align}
The mode probability vector evolves as
$p(t+\Delta t)=\exp(L\Delta t)p(t)$.
The associated semigroup $T_t=\exp(Lt)$ 
provides a finite-state approximation to the continuous Fokker--Planck equation
\begin{equation}
	\frac{\partial p(x,t)}{\partial t}
	=-\nabla\!\cdot(a(x)p)
	+\tfrac12\nabla^2\!:\!(\Sigma(x)p),
	\label{eq:fp}
\end{equation}
which governs the latent continuous dynamics prior to clustering.

The MFPT to a designated hazard set is computed
from the transient subgenerator $\tilde L$ (obtained by restricting $L$
to non-hazard modes) via $-\tilde L^{-1}\mathbf{1}$.

\subsection{Explicit Drift, Diffusion, and Generator Construction (Algorithm~\ref{alg:costate})}

The generator construction employed here is a moment-based, continuous-time
approximation consistent with diffusion limits under coarse temporal
sampling.

Define inter-cluster distances
\begin{align}
	l_{kl} &:= \frac{1}{n_k n_l}
	\sum_{i\in C_k}\sum_{j\in C_l} \|x_i-x_j\|_2, \\
	l_{kl}^{(q)} &:= \frac{1}{n_k n_l}
	\sum_{i\in C_k}\sum_{j\in C_l} \big|x_j^{(q)}-x_i^{(q)}\big|.
\end{align}
For transitions $C_k\to C_l$ ($k\neq l$),
\begin{align}
	a_{kl}^{(q)} &=
	\frac{1}{|T_{kl}|}
	\sum_{i\in T_{kl}}
	\frac{\big|x_{i+1}^{(q)}-x_i^{(q)}\big|}{t_{i+1}-t_i}, \\
	\bar a_k^{(q)} &=
	\sum_{l} \frac{a_{kl}^{(q)}}{|l_{kl}^{(q)}|}.
\end{align}
Define diffusion for $C_k\to C_l$ by
\begin{align}
	\sigma_{kl}^{pq}
	&=
	\frac{1}{|T_{kl}|-1}
	\sum_{i\in T_{kl}}
	\Big(
	\frac{\big|x_{i+1}^{(p)}-x_i^{(p)}\big|}{\sqrt{t_{i+1}-t_i}}
	- a_{kl}^{(p)}\sqrt{t_{i+1}-t_i}
	\Big)
	\nonumber\\
	&\quad\cdot
	\Big(
	\frac{\big|x_{i+1}^{(q)}-x_i^{(q)}\big|}{\sqrt{t_{i+1}-t_i}}
	- a_{kl}^{(q)}\sqrt{t_{i+1}-t_i}
	\Big),
\end{align}
with $\sigma_{kl}^{pq}=0$ if $|T_{kl}|\le 1$, and
$\bar\sigma_{kl}=\sum_{p,q}\sigma_{kl}^{pq}/(l_{kl}^{(p)}l_{kl}^{(q)})$,
set to zero if any denominator vanishes.

Construct
\begin{align}
	(L_1^{\mathrm{out}})_{lk} &= \bar a_{kl}\ (k\ne l),
	& (L_1^{\mathrm{out}})_{kk} &= -\sum_{l\ne k}(L_1^{\mathrm{out}})_{lk},\\
	(L_1^{\mathrm{in}})_{lk} &= -\bar a_{lk}\ (k\ne l),
	& (L_1^{\mathrm{in}})_{kk} &= -\sum_{l\ne k}(L_1^{\mathrm{in}})_{lk},\\
	(L_2)_{lk} &= \tfrac12(\bar\sigma_{kl}+\bar\sigma_{lk})\ (k\ne l),
	& (L_2)_{kk} &= -\sum_{l\ne k}(L_2)_{kl}.
\end{align}
Finally, set $L=L_1^{\mathrm{out}}+L_1^{\mathrm{in}}+L_2$ so that
$\mathbf{1}^\top L=0$ and
$p(t+\Delta t)=\exp(L\Delta t)p(t)$.

\paragraph*{Positivity and step size.}
Choose $\Delta t>0$ such that $\exp(L\Delta t)$ is elementwise
nonnegative. In practice, one enforces $L$ to be a valid generator
(off-diagonal nonnegative, row sums zero) and employs scaling-and-squaring
to compute $\exp(L\Delta t)$ when necessary.

\subsection{Generator Calibration and Validation}

To assess fidelity of the learned generator $L$, we employ:
\begin{enumerate}[leftmargin=3mm]
	\item \emph{One-step calibration error}:
	$\|\exp(L\Delta t)p_k-\hat p_{k+1}^{\mathrm{emp}}\|_1$, where $\hat p_{k+1}^{\mathrm{emp}}$ is the empirical mode histogram.
	\item \emph{Bootstrap confidence intervals} for entries $L_{kl}$ using
	resampled transition paths;
	\item \emph{Spectral stability}, monitoring the real parts of eigenvalues
	of $L$ to ensure contractive dynamics on the probability simplex.
\end{enumerate}

\section{Adaptive Weighting and Generator maximum likelihood estimates (MLE)}
\label{sec:AdaptiveWeighting}
\subsection{Adaptive $\sigma(x,t)$ from Information Metrics}

Observation weighting is selected using information-based criteria to
improve numerical conditioning and sensitivity of the co-state
correction. Let $S_t$ denote the empirical innovation covariance and
define
\begin{align}
	\sigma_t^{-2} \propto \alpha\, \mathcal{I}(x_t) + \beta I, \qquad
	\mathcal{I}(x_t) := H_t Q_t H_t^{\mathsf T},
	\label{eq:infosigma}
\end{align}
where $Q_t$ is a bounded positive semidefinite proxy for process
sensitivity (e.g., Jacobian-based), and $\alpha,\beta>0$ are tuning
parameters. This choice defines a stabilizing surrogate observation
model rather than a statistically optimal one.

Since $\sigma_t$ is computed from the current state and recent measurement history and is uniformly bounded, 
the adaptive weighting remains causal and numerically well behaved over finite horizons.

\subsection{Continuous-Time MLE}
The moment-based generator constructed in Section~VI may be further
refined via maximum likelihood estimation when sufficient transition
data are available. Let $N_{kl}(T)$ denote the number of observed jumps
from mode $k$ to $l$ over $[0,T]$, and let $T_k$ denote the total dwell
time in mode $k$. For a time-homogeneous generator, the log-likelihood is
\begin{align}
	\mathcal{L}(L)
	= \sum_{k\neq l} N_{kl}(T)\log L_{kl}
	- \sum_k L_{kk} T_k,
	\qquad
	L_{kk}=-\sum_{l\ne k}L_{kl}.
\end{align}
The maximization is performed under $L_{kl}\ge 0$ and
$\mathbf{1}^\top L=0$. Asymptotic covariance follows from the observed
information matrix, and MFPT uncertainty may be estimated via bootstrap
resampling of transition paths.

\section{Co-State Probability Correction}
\label{sec:Co-StateProbability}

For completeness within the closed-loop architecture, we briefly restate
the co-state probability correction used to enhance regime sensitivity.
Let the reachable co-state values generated by \eqref{eq:lamreg} be
clustered into a finite number of non-empty clusters with centroids
$\bar{\lambda}_k$. Given an observed co-state increment
$\Delta\hat{\lambda}$ over $\Delta t$, we apply the exponential
weighting correction based on the Radon-Nikodym derivative for changing the probability measure to that of the observed process and the particle representation of the Kallinapur-Striebel formula 
\cite{kallianpur1968,xiong2008}
\begin{equation}
	\hat{p}_k
	=
	\frac{
		\exp\!\left(
		\bar{\lambda}_k^{\mathsf T}\Delta\hat{\lambda}
		-\tfrac12\norm{\bar{\lambda}_k}^2\Delta t
		\right)
		p_k
	}{
		\sum_j
		\exp\!\left(
		\bar{\lambda}_j^{\mathsf T}\Delta\hat{\lambda}
		-\tfrac12\norm{\bar{\lambda}_j}^2\Delta t
		\right)
		p_j
	}.
	\label{eq:post}
\end{equation}
This correction augments the generator-driven probability evolution and
is applied selectively to enhance sensitivity under coarse or irregular
sampling. It should be viewed as a heuristic exponential weighting
motivated by Bayesian structure, rather than as an exact filtering
update. In practice, log-sum-exp evaluation and tempering thresholds are
used to ensure numerical stability.

\section{Detectability and Risk Awareness}

\subsection{Detectability for Consistency Convergence}

\begin{assumption}[Local Detectability]
	Along the trajectory, the linearization of $(f,H)$ is locally
	detectable in the standard control-theoretic sense, and the observation
	weighting $\sigma_t$ in \eqref{eq:wtobs} is bounded and computed causally
	from the current state and recent measurements.
\end{assumption}

\begin{proposition}[Consistency Integrability]
	Under Assumption~1 and bounded $H_t$, the consistency error functional
	$V(t)$ defined in Appendix~A is integrable, and co-state corrections occur
	only intermittently outside a compact neighborhood of the
	measurement-consistent manifold over finite horizons.
\end{proposition}

\subsection{Risk Awareness with Co-State and MFPT}
To illustrate how the proposed estimation and risk signals may be used
for decision-making, we embed them in a receding-horizon model predictive control (MPC) formulation.
Define the horizon-$N$ cost
\begin{equation}
	J
	=
	\sum_{k=0}^{N-1}
	\ell(x_k,u_k)
	+
	\gamma\|\lambda_k\|_1
	+
	\rho\,\mathrm{MFPT}(p_k)^{-1},
\end{equation}
subject to the probabilistic risk constraint
$p_k^\top W p_k \le R_{\max}$. The inverse MFPT term penalizes proximity to hazard states, while $\|\lambda_k\|_1$ discourages sustained geometric inconsistency.

State prediction uses the projected dynamics, and mode probabilities evolve according to $p_{k+1}=\exp(L\Delta t)p_k$ with optional correction via \eqref{eq:post}. The MPC is included to demonstrate how the proposed risk signals may be incorporated into a closed-loop architecture; no
claims of optimal control performance are made. The control layer is modular and does not affect the validity of the underlying estimation and risk prediction framework.

\begin{figure}[!ht]
	\centering
	\begin{tikzpicture}[node distance=10mm and 10mm]
		\node[draw,rounded corners,align=center,inner sep=3pt] (proj) {Projected\\ Dynamics};
		\node[draw,rounded corners,align=center,inner sep=3pt, right=17mm of proj] (gen) {Generator\\ $L$ Update};
		\node[draw,rounded corners,align=center,inner sep=3pt, below=9mm of gen] (corr) {Co-state\\ Correction};
		\node[draw,rounded corners,align=center,inner sep=3pt, right=10mm of gen] (mpc) {Risk-aware\\ MPC};
		\draw[->] (proj) -- node[above,sloped,near end]{\small $x_{k+1}$} (gen);
		\draw[->] (gen) -- node[right,near end]{\small $p_{k+1}$} (corr);
		\draw[->] (corr.east) -| ++(4mm,0) |- node[pos=0.25,above,sloped]{\small $\hat p_{k+1}$} (mpc.west);
		\draw[->] (mpc.west) |- ++(0mm,10mm) -| (proj.north);
		\node[left=4mm of proj] (y) {$dy^{\mathrm{obs}}$};
		\draw[->] (y) -- (proj);
	\end{tikzpicture}
	\caption{Risk-aware closed-loop architecture integrating projected dynamics, stochastic generator update, co-state correction, and model predictive control.}

\end{figure}
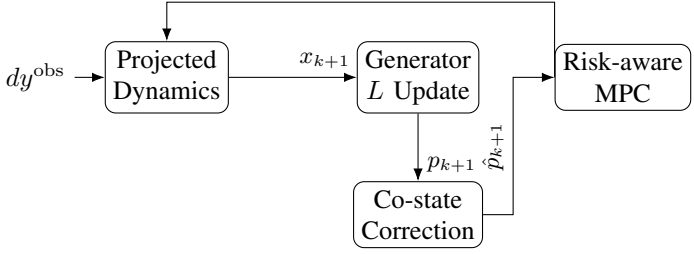

\section{Algorithm and Implementation}
\label{sec:Algorithm}
We summarize the proposed co-state fusion
pipeline as an executable algorithm suitable for real-time and
asynchronous onboard implementation.

\begin{algorithm}[!t]
	\caption{Co-State Based Data Fusion}
	\label{alg:costate}
	\begin{algorithmic}[1]
		\State \textbf{Input:} Measurement stream $\{(t_i,y_i)\}$, models $f,h$, noise model $R_t$
		\State \textbf{Initialize:} state $x_0$, co-state $\lambda_0=0$, mode probabilities $p_0$, generator $L_0$
		\State Compute $H_t$; update $(x_t,\lambda_t)$ via \eqref{eq:aug}--\eqref{eq:lamreg} (use regularized inverse if ill-conditioned)
		\State Stream $(x_t,\lambda_t)$; cluster online into modes $C_1,\dots,C_K$
		\State Estimate drifts $a_{kl}$ and diffusions $\Sigma_{kl}$; update generator $L$ (moment-based, with optional MLE refinement; Sec.~\ref{sec:Drift--Diffusion},~\ref{sec:AdaptiveWeighting}) using streamed state and co-state data
		\State Predict mode probabilities $p(t+\Delta t)=\exp(L\Delta t)p(t)$
		\State Correct probabilities via \eqref{eq:post} using observed co-state increment $\Delta\hat{\lambda}$
		\State \textbf{Asynchronous:} inject measurement pulses at sensor times; handle out-of-sequence measurements via bounded retroactive co-state updates over a finite buffer
		\State \textbf{Adaptive:} update observation weighting $\sigma_t$ via \eqref{eq:infosigma}
		\State \textbf{Control (optional):} using updated $(x_t,p_t)$, solve risk-aware MPC with constraint $p^\top W p\le R_{\max}$
		\State \textbf{Repeat online}
	\end{algorithmic}
\end{algorithm}

\begin{figure}[t] 
	\centering 
	\includegraphics[width=\linewidth]{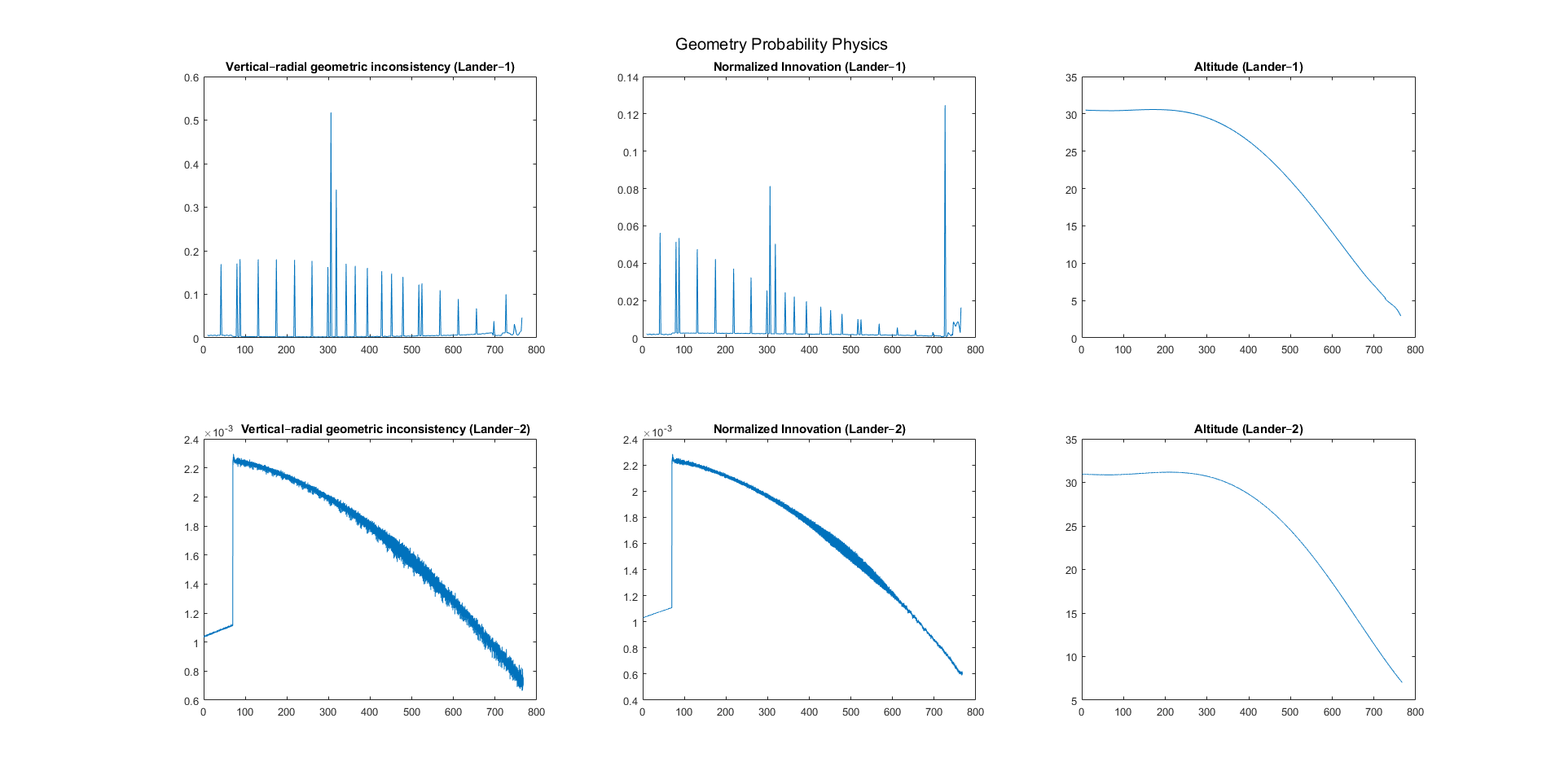} 
	\caption{Time histories of geometric inconsistency, normalized innovation, hazard probability, and physical descent states for Lander--1 and Lander--2.}

	\label{fig:geometry} 
\end{figure} 

\begin{figure}[t] 
	\centering 
	\includegraphics[width=\linewidth]{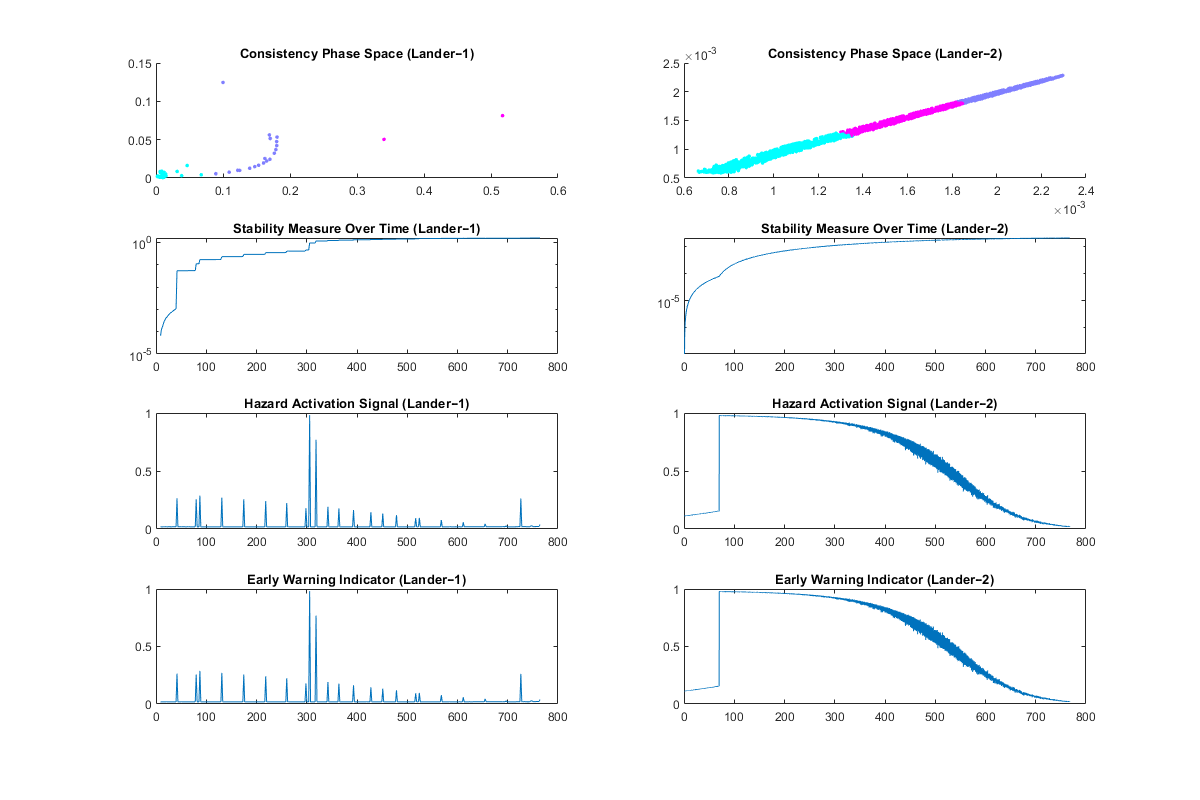} 
	\caption{Learned hazard probability, Lyapunov geometric stress, and early-warning indicators for both descents.} \label{fig:hazard} 
	
\end{figure} 

\begin{figure}[t] 
	\centering 
	\includegraphics[width=\linewidth]{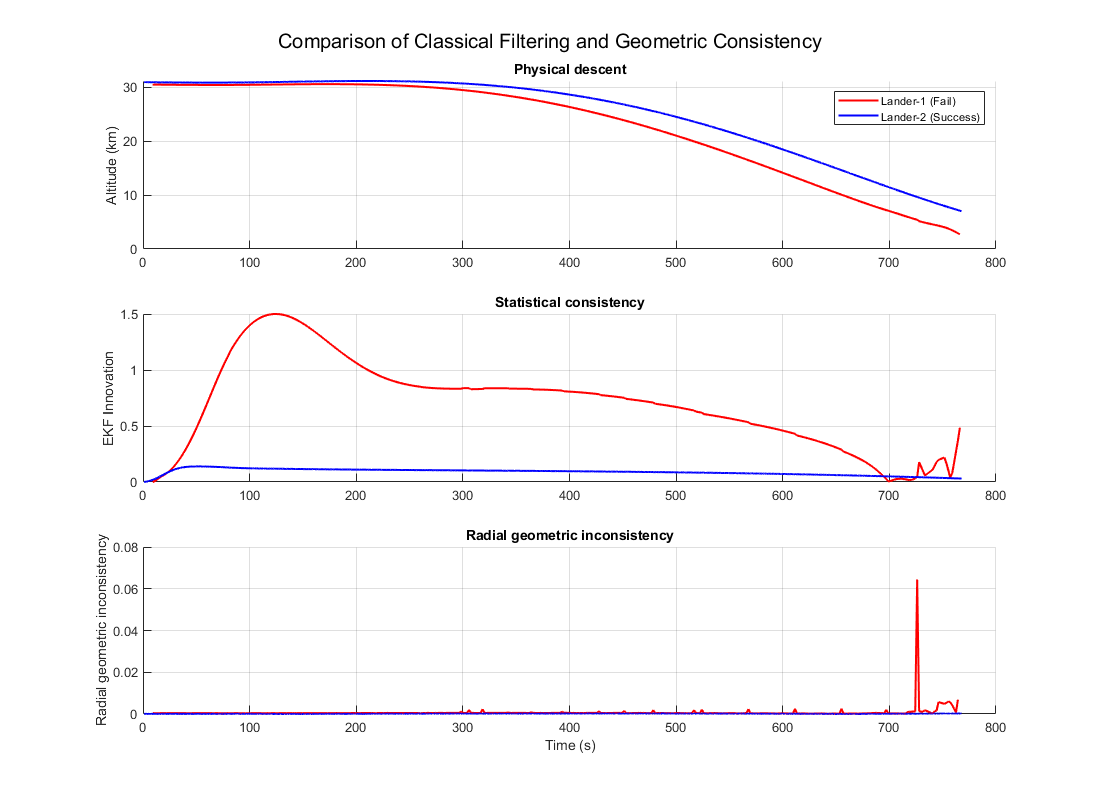} 
	
	\caption{Comparison of EKF innovation and co-state-based consistency signals for the failed descent (Lander--1). The co-state and whitened innovation rise earlier than the EKF innovation.}

	\label{fig:ekf-costate} 
\end{figure} 


\section{Results and Discussion: Lunar Lander Application}
\label{sec:results-lander}

This section evaluates the proposed co-state based fusion architecture on two real lunar powered-descent telemetry datasets\footnote{The telemetry logs correspond to Indian lunar landing missions and are anonymized for scientific neutrality and confidentiality.}. The descent spans approximately $30$~km altitude to touchdown and includes multiple operational phases with distinct navigation, guidance, and control (NGC) characteristics. \emph{Lander--1} corresponds to a failed descent whose behavior is consistent with an internal modeling or software inconsistency within the guidance–propulsion chain, while \emph{Lander--2} corresponds to a successful autonomous landing.

The purpose of this section is not to assess state-estimation accuracy, but to demonstrate how the geometric and stochastic structures developed in Sections~\ref{sec:sdamodel}--\ref{sec:AdaptiveWeighting} manifest in real descent telemetry and expose latent internal inconsistencies well before any physical instability, loss of control authority, or estimator divergence becomes visible.

Specifically, we demonstrate that:
(i) the co-state differential--algebraic equation exposes instantaneous violations of the assumed measurement--dynamics geometry;
(ii) repeated violations accumulate coherently across guidance update cycles and induce slow-time drift toward hazardous regimes in the learned stochastic generator;
(iii) this drift produces early, interpretable risk signals without labeled failures or tuned thresholds; and
(iv) the resulting diagnostics isolate internal modeling inconsistency rather than physical or navigational failure.

Throughout this section, the \emph{co-state magnitude} denotes the norm of the algebraic multiplier enforcing measurement--dynamics consistency, the \emph{normalized innovation} denotes the whitened measurement residual, and the \emph{hazard probability} denotes the probability of occupying the hazard regime under the learned continuous-time generator.

\subsection{Mission Phasing and Signal Construction}

Powered descent proceeds through three operational regimes:
(i) high-altitude inertial descent ($30$--$7$~km), dominated by propagation with limited measurement updates;
(ii) mid-altitude transition ($7$--$2$~km), where absolute navigation sensors activate and guidance begins incorporating measurements; and
(iii) terminal descent ($<$2~km), where guidance and control operate in tightly coupled closed-loop form.

The translational state is
\[
x = [x,y,z,v_x,v_y,v_z]^\top,
\]
and the measurement vector is chosen as
\[
y^{\mathrm{obs}} = [z,\ \|p\|,\ v_z]^\top, \qquad p=[x,y,z]^\top,
\]
with analytic Jacobian $H_t=\partial h/\partial x(x_t,t)$.

The nominal dynamics underlying both navigation and guidance assume the small-noise kinematic model $\dot p=v,\ \dot v=0$, implicitly assuming that commanded thrust tracks desired acceleration and that unmodeled forces and actuation nonlinearities are negligible.

The instantaneous geometric correction is given by
\[
\lambda_t\,dt = (H_tH_t^\top)^{-1}\Sigma_t^{-1} \big(dy_t^{\mathrm{obs}} - \eta_t\,dt\big),
\]
where $\eta_t$ is the model-predicted measurement increment and $\Sigma_t$ whitens innovations. Under nominal consistency, the normalized innovation
\[
z_t=\|\Sigma_t^{-1/2}(dy_t^{\mathrm{obs}}-\eta_t\,dt)\|_2
\]
behaves approximately as a $\chi_3$ variable.

The feature vector $\phi_t=(\lambda_t,\|\lambda_t\|,z_t)$ is clustered into \emph{Nominal}, \emph{Corrective}, and \emph{Hazard} regimes, yielding a continuous-time generator $L$ that captures slow-time regime evolution.

Note that two related but distinct geometric inconsistency signals are used in this study. The EKF comparison uses a radial-geometry co-state enforcing consistency of altitude, range, and radial velocity, while the full fusion architecture uses a vertical–radial co-state enforcing consistency of vertical position, range, and vertical velocity. Both arise from the same projection principle, but operate on different measurement manifolds and should be interpreted as complementary geometric residuals rather than identical quantities.

\subsection{Deterministic Geometry: Internal Consistency of NGC}

Figure~\ref{fig:geometry} shows the co-state magnitude and normalized innovation for both descents.

\paragraph*{Lander--1 (failed).}
From high altitude onward, the co-state magnitude remains persistently nonzero, indicating that the assumed kinematic and thrust-response model cannot simultaneously satisfy the measurement geometry. As navigation updates become active in the $7$--$2$~km regime, this inconsistency grows and exhibits intermittent sharp excursions.

These excursions do not correspond to any physical instability: altitude and velocity evolve smoothly, thrust remains effective, and control authority is maintained. Instead, the inconsistency reflects that the internally propagated motion model used by guidance and navigation does not match the effective dynamics implied by the measurements. In practice, this is consistent with a software-level model mismatch, such as incorrect thrust mapping, discretization error, saturation modeling error, or timing misalignment between guidance and propulsion execution.

Thus, the co-state reveals a structural inconsistency internal to the NGC system with respect to the assumed measurement–dynamics geometry, as defined by the active navigation observables.

\paragraph*{Lander--2 (successful).}
In contrast, the co-state remains near zero throughout descent. Minor excursions following navigation updates collapse rapidly, indicating that measurement updates, guidance computation, and thrust execution remain mutually compatible. The system remains close to a single consistent measurement--dynamics manifold, reflecting coherent operation of navigation, guidance, and control.

\subsection{Stochastic Regimes: Accumulation of Structural Stress}

Figure~\ref{fig:hazard} shows the regime probabilities inferred from the generator.

For Lander--2, brief transitions into corrective regimes occur after navigation updates but rapidly return to nominal, indicating that corrections are absorbed by the closed-loop system without structural stress.

For Lander--1, however, repeated transitions into hazard regimes occur throughout mid- and terminal descent. Each transition corresponds to a period in which the deterministic geometry is violated; although each violation is transient, their repeated occurrence induces slow-time drift toward hazardous regimes in the generator. Importantly, this drift arises without any abrupt physical anomaly, reflecting accumulation of internal modeling stress rather than physical instability.

\subsection{Comparison with Classical EKF Diagnostics}

Figure~\ref{fig:ekf-costate} compares physical descent, EKF innovation, and geometric inconsistency.

The EKF remains numerically stable in both cases. For Lander--1, however, its innovation statistics become inconsistent only late in the descent, after substantial internal inconsistency has already accumulated. In contrast, the co-state responds immediately to geometric violations, providing earlier and more interpretable diagnostic information.

This highlights the distinction between statistical consistency (as monitored by EKF residuals) and geometric consistency (as enforced by the co-state). The former reacts to accumulated errors, while the latter exposes structural violations at their point of origin.

\subsection{Interpretation and Implications}

Taken together, the results indicate that:
(i) the failure is not caused by physical instability, loss of thrust, or navigation error;
(ii) it originates from an internal modeling or software inconsistency in the guidance–propulsion chain;
(iii) this inconsistency manifests first as a geometric violation, not as a statistical anomaly or physical divergence; and
(iv) the proposed architecture detects and aggregates this violation into a probabilistic risk signal well before terminal failure.

Thus, the deterministic–stochastic co-state fusion framework provides a principled mechanism for detecting internal system-level inconsistencies that are invisible to traditional state-estimation diagnostics but critical for safety and autonomy in planetary landing systems.

\section{Discussion and Future Work}
\label{sec:DiscussionFuturework}

This work introduced a co-state based fusion architecture that enforces instantaneous measurement--dynamics consistency via a differential--algebraic geometric projection, while capturing slow-time hazard evolution through data-driven stochastic generator learning.

At the deterministic level, the algebraic co-state acts as an instantaneous Lagrange multiplier that projects the nominal state evolution onto the measurement-consistent tangent manifold. This yields a stochastic differential--algebraic system whose trajectories remain locally compatible with the sensing geometry even in the presence of modeling bias, partial observability, or ill-conditioned measurement Jacobians. In the small-noise regime, the formulation reduces to the deterministic limit of innovation-based nonlinear filtering, while preserving an explicit geometric interpretation that is not accessible through purely statistical residuals.

At the stochastic level, the learned continuous-time generator provides a finite-dimensional approximation of latent regime dynamics induced by accumulated geometric inconsistency and model mismatch. This representation enables intrinsic probabilistic hazard forecasting through mode occupancy probabilities, spectral properties of the generator, and derived quantities such as mean first-passage time (MFPT). The empirical results demonstrate that geometric inconsistency, stochastic drift, and risk measures rise coherently prior to failure, yielding operationally meaningful early-warning capability on real lunar descent telemetry.

Importantly, the deterministic projection and stochastic update operations are low-dimensional and constant-time per step, making the framework computationally comparable to or cheaper than a standard EKF update and therefore feasible for real-time onboard execution.

Several directions remain open for future work.

\paragraph*{Risk-aware guidance and control integration.}
The co-state magnitude, hazard probability, and MFPT can be explicitly incorporated into guidance and control laws, enabling adaptive risk-aware behavior in which control authority, descent aggressiveness, or replanning triggers are modulated based on predicted hazard proximity and geometric stress. Embedding the framework within receding-horizon model predictive control formulations would enable closed-loop mitigation of emerging hazards rather than purely diagnostic monitoring.

\paragraph*{Real-time and embedded implementation.}
The algebraic structure of the co-state update and the low-order nature of the learned generator make the framework amenable to real-time onboard implementation. Cholesky-based solvers for the regularized projection, sparse or structured generator updates, and low-rank or Krylov approximations for matrix exponentials suggest feasible deployment on flight processors, FPGAs, or system-on-chip platforms under stringent computational and power constraints.

\paragraph*{Multi-agent and cooperative descent.}
The framework naturally extends to multi-agent and multi-vehicle scenarios, in which synchronized co-states can encode shared geometric constraints, mutual observability, and collective hazard awareness. This opens the possibility of cooperative landing, distributed fault detection, and shared risk assessment in formation-flying or swarm-based planetary exploration missions.

\paragraph*{Perception and environment coupling.}
Integration with vision-based terrain-relative navigation, hazard detection, and environmental perception would enable coupling between internal geometric consistency and external environmental risk. This would allow the system to reason jointly about internal model mismatch and externally induced hazards, supporting end-to-end risk-aware autonomy for unstructured planetary environments.

\paragraph*{Validation, robustness, and certification.}
Systematic validation across diverse missions and failure modes, including ROC-based detection analysis, uncertainty quantification of generator parameters and MFPT, and robustness studies under sensor dropout, delays, and adversarial disturbances, are essential for operational readiness. Formal connections to detectability, observability, and fault isolation theory would further strengthen theoretical guarantees and support eventual certification.

Taken together, these directions position the proposed co-state fusion framework not only as a diagnostic mechanism, but as a foundation for future geometry-consistent, probabilistically grounded, and risk-aware autonomous navigation and control systems.

\appendices

\section{Proof of Lyapunov Stability}

Define the consistency error functional
\begin{equation}
	V(t)
	:= \tfrac12
	(dx_t - f_t dt)^\top
	(H_t H_t^\top + \varepsilon I)^{-1}
	(dx_t - f_t dt),
	\label{eq:lyap}
\end{equation}
which measures deviation from the measurement-consistent tangent manifold.
The functional is interpreted in the It\^{o} sense, accounting for quadratic
variation terms induced by the stochastic observation process.

Under boundedness and regularity assumptions on $H_t$, $\sigma_t$, and $f$,
differentiation of $V(t)$ along trajectories of
\eqref{eq:aug}--\eqref{eq:lamreg} yields
\[
\dot V(t)
=
-\lambda_t^{\mathsf T}
(H_t H_t^{\mathsf T})
\lambda_t
+ \text{martingale terms}.
\]
Taking expectations and using the positive semidefiniteness of $H_t H_t^\top$,
the martingale terms have zero expectation, yielding
\[
\mathbb{E}[\dot V(t)] \le 0.
\]

Thus, the consistency error is non-increasing in expectation and integrable
over finite horizons. Existence and uniqueness of solutions for the resulting
stochastic differential--algebraic system follow from standard index-1 DAE
theory \cite{ascher1998dae,kunkel2006dae}.

\section{Regularized Projector}

Let $P_\perp = H_t^\top (H_t H_t^\top)^{-1} H_t$ denote the orthogonal projector
onto $\Range(H_t^\top)$. The regularized projector
$P_\perp^\varepsilon := H_t^\top (H_t H_t^\top + \varepsilon I)^{-1} H_t$
satisfies
\[
\|P_\perp^\varepsilon - P_\perp\|_2 = \mathcal{O}(\varepsilon),
\qquad \varepsilon \to 0,
\]
provided $H_t$ remains bounded. This ensures numerical stability under rank
deficiency while preserving the geometric meaning of the projection.

\section{Detectability Sketch}

Under local detectability of the linearized pair $(f,H)$ along the trajectory,
the innovation remains informative along unstable directions. Bounded and
causally computed observation weighting $\sigma_t$ ensures numerical stability.

If the consistency error were to persist outside a compact neighborhood of the
measurement-consistent manifold, it would contradict detectability via
persistent excitation of the innovation. Consequently, the error functional
$V(t)$ is integrable, and large co-state corrections occur only intermittently
over finite horizons.

\section{Implementation Notes}

In practice, $(H_t H_t^\top)^{-1}$ is computed via Cholesky factorization.
When ill-conditioning is detected, the regularized inverse
$(H_t H_t^\top + \varepsilon I)^{-1}$ is used, with $\varepsilon$ selected
adaptively based on spectral thresholds.

The continuous-time generator $L$ is estimated via constrained maximum
likelihood with $L_{kl}\ge 0$ for $k\neq l$ and $\mathbf{1}^\top L=0$.
Matrix exponentials are evaluated using scaling-and-squaring, and Bayesian
updates are implemented in log-domain for numerical stability.

For embedded implementations, explicit matrix inversion is avoided, sparsity
patterns are precomputed, and generator propagation may be approximated using
low-order rational schemes.

\bibliographystyle{IEEEtran}
\bibliography{references_v4}

@article{kalman1960,
  author = {Kalman, R. E.},
  title = {A New Approach to Linear Filtering and Prediction Problems},
  journal = {ASME Trans., J. Basic Eng.},
  year = {1960},
  volume = {82},
  number = {1},
  pages = {35--45}
}

@article{kalmanbucy1961,
  author = {Kalman, R. E. and Bucy, R. S.},
  title = {New Results in Linear Filtering and Prediction Theory},
  journal = {ASME Trans., J. Basic Eng.},
  year = {1961},
  volume = {83},
  number = {1},
  pages = {95--108}
}

@book{doucet2001,
  editor = {Doucet, Arnaud and de Freitas, Nando and Gordon, Neil},
  title = {Sequential Monte Carlo Methods in Practice},
  publisher = {Springer},
  year = {2001}
}

@book{ascher1998dae,
  author = {Ascher, Uri M. and Petzold, Linda R.},
  title = {Computer Methods for Ordinary Differential Equations and Differential-Algebraic Equations},
  publisher = {SIAM},
  year = {1998}
}

@book{kunkel2006dae,
  author = {Kunkel, Peter and Mehrmann, Volker},
  title = {Differential-Algebraic Equations: Analysis and Numerical Solution},
  publisher = {EMS Press},
  year = {2006}
}

@article{kallianpur1968,
  author = {Kallianpur, G. and Striebel, C.},
  title = {Estimation of Stochastic Systems: Arbitrary System Process with Additive Noise Observation Errors},
  journal = {Ann. Math. Stat.},
  year = {1968},
  volume = {39},
  number = {3},
  pages = {785--801}
}

@article{GaoZhou2018_TerrainHazard,
  author  = {Gao, Ai and Zhou, Shibo},
  title   = {Feature Density Based Terrain Hazard Detection for Planetary Landing},
  journal = {IEEE Transactions on Aerospace and Electronic Systems},
  year    = {2018},
  volume  = {54},
  number  = {5},
  pages   = {2411--2420},
  doi     = {10.1109/TAES.2018.2818578}
}

@article{Fourlas2021_FaultDiagnosisFTC,
  author  = {G. K. Fourlas},
  title   = {A Survey on Fault Diagnosis and Fault-Tolerant Control Methods for Unmanned Aerial Vehicles},
  journal = {Aerospace},
  year    = {2021},
  volume  = {9},
  number  = {9},
  pages   = {197},
  doi     = {10.3390/aerospace9090197}
}

@article{WangPengJiang2020_RT_FaultDetection,
  author  = {B. Wang and X. Peng and M. Jiang and D. Liu},
  title   = {Real-Time Fault Detection for UAV Based on Model Acceleration Engine},
  journal = {IEEE Transactions on Instrumentation and Measurement},
  year    = {2020},
  volume  = {69},
  pages   = {9505--9516}
}

@article{JungBang2021_FaultTolerantMPC,
  author    = {W. Jung and H. Bang},
  title     = {Fault and Failure Tolerant Model Predictive Control of Quadrotor UAV},
  journal   = {International Journal of Aeronautical and Space Sciences},
  year      = {2021},
  volume    = {22},
  number    = {3},
  pages     = {663--675}
}

@conference{Calkins2022_RobustTrajectoryPDL,
  author    = {G. E. Calkins and D. C. Woffinden and Z. R. Putnam},
  title     = {Robust Trajectory Optimization for Guided Powered Descent and Landing},
  booktitle = {AAS/AIAA Astrodynamics Specialists Conference},
  year      = {2022},
  note      = {Paper AAS 22-660, NASA NTRS 20220010431}
}

@book{xiong2008,
  title={An introduction to stochastic filtering theory},
  author={Xiong, Jie},
  volume={18},
  year={2008},
  publisher={OUP Oxford}
}

\end{document}